\pdfoutput=1

\documentclass[sigconf,9pt,nonacm]{acmart}

\setcopyright{acmcopyright}
\copyrightyear{2018}
\acmYear{2018}
\acmDOI{XXXXXXX.XXXXXXX}

\acmConference[Conference acronym 'XX]{Make sure to enter the correct
  conference title from your rights confirmation emai}{June 03--05,
  2018}{Woodstock, NY}
\acmPrice{15.00}
\acmISBN{978-1-4503-XXXX-X/18/06}

\usepackage{xifthen}
\usepackage{xspace}
\usepackage{tikz-cd}

\usepackage{xurl}
\usepackage{setouterhbox}
\usepackage{etoolbox}

\makeatletter
\patchcmd\deferred@thm@head{\sbox\@labels{\normalfont#1}}{\setouterhbox\@labels \normalfont#1\endsetouterhbox}{}{\fail}
\makeatother

\usetikzlibrary{decorations.markings}
\tikzset{mid vert/.style={/utils/exec=\tikzset{every node/.append style={outer sep=0ex}},
		postaction=decorate,decoration={markings,
			mark=at position 0.5 with {\draw[-] (0,#1) -- (0,-#1);}}},
	mid vert/.default=0.75ex}

\newcommand{\longhash}{69ecbf1322dfb67cdb5a5dbbe864b22695d2e60b}
\newcommand{\shorthash}{69ecbf1}

\newcommand{\coqdocbasebaseurl}{https://anonymous123.gitlab.io/UniMath-doc/\shorthash}

\newcommand{\coqdocbaseurl}{\coqdocbasebaseurl /UniMath.}
\newcommand{\urlhash}{\#}
\newcommand{\coqdocurl}[2]{\coqdocbaseurl #1.html\urlhash #2}
\newcommand{\nolinkcoqident}[1]{\normalfont\color{blue}{\nolinkurl{#1}}} %
\makeatletter
\newcommand{\coqident}{\begingroup\@makeother\#\@coqident}
\newcommand{\@coqident}[3][]{%
  \ifthenelse{\isempty{#2}}%
  {\nolinkcoqident{#3}}%
  {\ifthenelse{\isempty{#1}}%
  {\href{\coqdocurl{#2}{#3}}{\nolinkcoqident{#3}}}%
  {\href{\coqdocurl{#2}{#3}}{\nolinkcoqident{#1}}}}%
\endgroup}
\newcommand{\coqfile}[2]{%
  \ifthenelse{\isempty{#1}}%
  {\href{\coqdocbaseurl #2.html}{\nolinkcoqident{#2.v}}}%
  {\href{\coqdocbaseurl #1.#2.html}{\nolinkcoqident{#2.v}}}}
\makeatother

\AtEndPreamble{%
\theoremstyle{acmdefinition}
\newtheorem{constrInternal}[theorem]{Construction}%
\newtheorem{remark}[theorem]{Remark}
}

\usepackage[draft]{optional}
\newcommand{\plan}[1]{}
\newcommand{\BA}[1]{}
\newcommand{\PN}[1]{}
\newcommand{\NR}[1]{}
\newcommand{\NW}[1]{}
\opt{draft}{
\renewcommand{\plan}[1]{\textcolor{blue}{TODO: #1}\PackageWarning{TODO}{TODO: #1}}
\renewcommand{\BA}[1]{\textcolor{orange}{BA: #1}\PackageWarning{TODO}{TODO: #1}}
\renewcommand{\PN}[1]{\textcolor{purple}{PN: #1}\PackageWarning{TODO}{TODO: #1}}
\renewcommand{\NR}[1]{\textcolor{lime}{NR: #1}\PackageWarning{TODO}{TODO: #1}}
\renewcommand{\NW}[1]{\textcolor{magenta}{NW: #1}\PackageWarning{TODO}{TODO: #1}}
}

\newcommand{\teletype}[1]{\ensuremath{\mathtt{#1}}}
\newcommand{\systemname}[1]{\teletype{\color{darkgray}#1}\xspace}
\newcommand{\UniMath}{\systemname{UniMath}}

\newcounter{saveenumi}

\usepackage[capitalize]{cleveref}

\newcommand{\constfont}[1]{\ensuremath{\mathsf{#1}}}
\newcommand{\cat}[1]{\ensuremath{\constfont{#1}}\xspace}
\newcommand{\CC}[0]{\cat{C}} %
\newcommand{\CD}[0]{\cat{D}} %
\newcommand{\CV}[0]{\cat{V}} %
\newcommand{\SET}[0]{\cat{Set}} %
\newcommand{\id}[1]{\cat{id}_{#1}} %
\newcommand{\nattrans}[2]{#1 \Rightarrow #2} %
\newcommand{\iso}[2]{#1 \cong #2} %
\newcommand{\op}[1]{#1^{\cat{op}}} %
\newcommand{\idtoiso}{\cat{idtoiso}} %

\newcommand{\CB}[0]{\cat{B}} %
\newcommand{\cell}[2]{#1 \Rightarrow #2} %
\newcommand{\co}[1]{#1^{\cat{co}}} %
\newcommand{\runitor}{\rho}

\newcommand{\linvunitor}{\ell^{-1}}
\newcommand{\adjunction}[2]{#1 \dashv #2}

\newcommand{\VB}[0]{\cat{B}} %
\newcommand{\HorB}[1][]{\cat{H}_{#1}} %
\newcommand{\VerB}[1][]{\cat{V}_{#1}} %
\newcommand{\HorM}[2]{#1 \rightarrow #2} %
\newcommand{\VerM}[2]{#1 \rightarrow #2} %
\newcommand{\HorC}[2]{#1 \Rightarrow #2} %
\newcommand{\VerC}[2]{#1 \Rightarrow #2} %
\newcommand{\lwhisker}[2]{#1 \mathop{\vartriangleleft} #2} %
\newcommand{\rwhisker}[2]{#1 \mathop{\vartriangleright} #2} %
\newcommand{\dwhisker}[2]{#1 \mathop{\triangledown} #2} %
\newcommand{\uwhisker}[2]{#1 \mathop{\vartriangle} #2} %

\newcommand{\CellToSquareH}{\cat{CellToSq}_{\cat{H}}} %
\newcommand{\CellToSquareV}{\cat{CellToSq}_{\cat{V}}} %

\newcommand{\idtogreg}[2]{\cat{IdToGregEq}_{#1, #2}} %

\newcommand{\SquareV}[1]{\cat{Sq}(#1)}
\newcommand{\LAdj}[1]{\cat{LAdj}(#1)}
\newcommand{\DoubleCatToVerity}[1]{\overline{#1}}
\newcommand{\ProfV}{\cat{Prof}}
\newcommand{\VProfV}[1]{\cat{Prof}_{#1}}
\newcommand{\Coreflect}[1]{\cat{Coreflect}(#1)}
\newcommand{\Polys}[1]{\cat{Poly}_{#1}}

\newcommand{\DC}{\cat{C}}
\newcommand{\VertCat}[1]{\cat{Ver}(#1)}
\newcommand{\HorCat}[1]{\cat{Hor}(#1)}
\newcommand{\SquareS}[1]{\cat{Sq}(#1)}
\newcommand{\SpanS}[1]{\cat{Span}(#1)}
\newcommand{\SquareU}[1]{\cat{Sq}(#1)}
\newcommand{\ProfDC}{\cat{Prof}}
\newcommand{\BicatToDoubleCat}[1]{\widehat{#1}}

\newcommand{\CompanionUnit}{\eta}
\newcommand{\CompanionCounit}{\varepsilon}

\newcommand{\profunctor}[2]{#1 \rightarrow #2}
\newcommand{\reprL}[1]{\cat{rep}_{\ell}(#1)}

\newcommand{\midhor}{
	\hspace{-0.5em}
  \begin{tikzcd}[sep=0.2in]
		 \arrow[r, mid vert] \pgfmatrixnextcell \
	\end{tikzcd}
  \hspace{-.7em}
}
\newcommand{\ver}[2]{{#1} \rightarrow {#2}} %
\newcommand{\hor}[2]{{#1} \midhor {#2}} %
\newcommand{\horid}[1]{\cat{id}_{#1}} %
\newcommand{\horcomp}[2]{#1 \odot #2} %
\renewcommand{\rightarrow}{
  \hspace{-0.5em}
  \begin{tikzcd}[sep=0.2in]
		 \arrow[r] \pgfmatrixnextcell \
	\end{tikzcd}
  \hspace{-.7em}
}

\newcommand{\Dsquare}[4]{\left( #1 \, \, \, ^{#3}_{#4} \, \, \, #2 \right)} %
\newcommand{\sqvid}[1]{\cat{id}^{v}_{\cat{sq}}(#1)} %
\newcommand{\sqvcomp}[2]{#1 \cdot_{\cat{sq}} #2} %
\newcommand{\sqhid}[1]{\cat{id}^{h}_{\cat{sq}}(#1)} %
\newcommand{\sqhcomp}[2]{#1 \odot_{\cat{sq}} #2} %
\newcommand{\Dlunitor}[1]{\lambda_{#1}}
\newcommand{\Drunitor}[1]{\rho_{#1}}
\newcommand{\Dassociator}[3]{\alpha_{(#1, #2, #3)}}

\newcommand{\StructCospan}[1]{\cat{StructCospan}(#1)}

\newcommand{\SpanDiag}[5]{#1 \xleftarrow{#2} #3 \xrightarrow{#4} #5} %
\newcommand{\Cospan}[5]{#1 \xrightarrow{#2} #3 \xleftarrow{#4} #5} %

\newcommand{\fullsubbicat}[2]{} %

\newcommand{\bicatofcats}[0]{\cat{UnivCat}} %
\newcommand{\bicatofEcats}[1]{\cat{UnivCat}_{#1}} %

\newcommand{\doublecatvercat}[1]{} %
\newcommand{\doublecathor}[1]{} %
\newcommand{\doublecatlunitor}[1]{} %
\newcommand{\doublecatrunitor}[1]{} %
\newcommand{\doublecatassociator}[1]{} %

\begin{document}

\title{Insights From Univalent Foundations: A Case Study Using Double Categories}

\author{Nima Rasekh}
\authornote{Both authors contributed equally to this research.}
\email{rasekh@mpim-bonn.mpg.de}
\orcid{0000-0003-0766-2755}
\affiliation{%
  \institution{Max Planck Institute for Mathematics}
  \city{Bonn}
  \country{Germany}
}

\author{Niels van der Weide}
\authornotemark[1]
\email{nweide@cs.ru.nl}
\orcid{0000-0003-1146-4161}
\affiliation{%
  \institution{Radboud University}
  \city{Nijmegen}
  \country{The Netherlands}
}

\author{Benedikt Ahrens}
\email{B.P.Ahrens@tudelft.nl}
\orcid{0000-0002-6786-4538}
\affiliation{%
  \institution{Delft University of Technology}
  \country{The Netherlands}
}
\affiliation{%
  \institution{University of Birmingham}
  \country{United Kingdom}
}

\author{Paige Randall North}
\email{p.r.north@uu.nl}
\orcid{0000-0001-7876-0956}
\affiliation{%
  \institution{Utrecht University}
  \country{The Netherlands}
}

\renewcommand{\shortauthors}{Rasekh, Van der Weide, Ahrens, and North}

\begin{abstract}
Category theory unifies mathematical concepts, aiding comparisons across structures by incorporating objects and morphisms, which capture their interactions. It has influenced areas of computer science such as automata theory, functional programming, and semantics. Certain objects naturally exhibit two classes of morphisms, leading to the concept of a \emph{double category}, which has found applications in computing science (e.g., ornaments, profunctor optics, denotational semantics).

The emergence of diverse categorical structures motivated a unified framework for category theory. However, unlike other mathematical objects, classification of categorical structures faces challenges due to various relevant equivalences. This poses significant challenges when pursuing the formalization of categories and restricts the applicability of powerful techniques, such as transport along equivalences. This work contends that \emph{univalent foundations} offers a suitable framework for classifying different categorical structures based on desired notions of equivalences, and remedy the challenges when formalizing categories. The richer notion of equality in univalent foundations makes the equivalence of a categorical structure an inherent part of its structure.  

We concretely apply this analysis to double categorical structures. We characterize and formalize various definitions in Coq UniMath, including (pseudo) double categories and double bicategories, up to chosen equivalences. We also establish \emph{univalence principles}, making chosen equivalences part of the double categorical structure, analyzing strict double setcategories (invariant under isomorphisms), pseudo double setcategories (invariant under isomorphisms), univalent pseudo double categories (invariant under vertical equivalences) and univalent double bicategories (invariant under gregarious equivalences).

\end{abstract}



\received{20 February 2007}
\received[revised]{12 March 2009}
\received[accepted]{5 June 2009}

\maketitle

\section{Introduction}
\label{sec:introduction}

This work focuses on the fundamental connection between \emph{category theory}, \emph{univalent foundations} and the \emph{formalization of mathematics}.
In this introductory section we demonstrate that developing category theory with classical foundations limits its applicability in mathematics (\ref{subsec:math}) and computer science (\ref{subsec:cs}).
We then explain how univalent foundations, and particularly the \emph{univalence maxim} for category theory (\ref{subsec:univ}) that we introduce, overcomes these limitations.
We justify our assertion by applying the maxim to several categorical notions, including categories, $2$-categories and double categories (\ref{subsec:results}).

\subsection{Limitation of Classical Foundations for Category Theory in Mathematics}
\label{subsec:math}
A category can be seen as an inherently algebraic structure, generalizing monoids. As such, categories -- like groups and other algebraic objects -- come with a natural notion of sameness commonly known as \emph{isomorphisms of categories}. However, with the advent of category theory, Eilenberg and Mac Lane also introduced a relevant notion not commonly found in algebra -- namely that of natural transformation -- which results in a weaker notion of sameness: \emph{equivalences of categories}. This early dichotomy between two notions of sameness that naturally arise out of one definition permeates all of the category theory literature.  

In recent decades categories have been generalized in a variety of manners under the banner of \emph{higher category theory} or \emph{higher-dimensional category theory}. Examples of such structures include $2$-categories and more generally $n$-categories, but also double categories, and even $(\infty,n)$-categories. As one might expect, the ambiguity outlined in the previous paragraph is exacerbated in this situation and every one of these higher categorical structures naturally comes with a wide range of notions of sameness. For example, $2$-categories have at least $5$ notions of sameness ranging from \emph{isomorphisms} up to \emph{biequivalences}. This suggests the need to develop ways to resolve such ambiguities when working with categorical data. 

In classical mathematical foundations the objects by themselves do not give us any information about the desired equivalences. Hence, the first proposed solution consists of choosing the context appropriately. For example, we can consider categories as objects in a $1$-category of categories, explicitly ignoring natural transformations and hence invariant under isomorphisms, or as objects in a $2$-category, making it invariant under equivalences. While this approach is an improvement upon the previous state of affairs, it leaves us with certain challenges. Indeed, it is certainly the case that even if we study functors between isomorphism-invariant categories, we would often want to consider their natural transformations. Moreover, it also creates a communicative barrier where every object needs to constantly be explicitly presented in its context, lest ambiguities arise. Finally, we face a mismatch between the built-in equality of the foundations and the chosen equivalences. For example, we would generally not expect the set of objects of two equivalent categories to be equal (or even isomorphic).

\subsection{Limitation of Classical Foundations for Category Theory in Computer Science}
\label{subsec:cs}
\emph{Formalization of mathematics} is a pivotal meeting ground for computer science and mathematics, contributing towards progress in theoretical computer science and formal verification of sophisticated mathematics. Its importance has significantly grown, due to advances in interactive theorem provers and the intricate landscape of contemporary mathematics. The central role of category theory in modern mathematics necessitates an effective approach to the formalization of categorical concepts that is both computationally feasible and applicable across various formal systems.

Pursuing the formalization of a sophisticated mathematical object such as a category necessitates appropriate type theoretical and foundational assumptions. Indeed, unlike groups and other familiar algebraic objects, categories are most naturally described as dependent types, which by construction include source and target of a morphism and simplify the data of a composition, unitality and associativity. On the other hand, it is well known that working with extensional identity types results in not having decidable type checking \cite{hofmann1995extensional}. As a result, the formalization of category theory has primarily, however not exclusively \cite{bayer2023category}, advanced in proof assistants relying on dependent type theories with intensional identity types, such as Agda \cite{DBLP:conf/cpp/HuC21}, Coq \cite{rezk_completion}, or Lean \cite{zbMATH07606352}.

Unfortunately, a straightforward approach towards a formalization of category theory in such settings faces numerous challenges, primarily stemming from the multifaceted nature of equivalences inherent to categories. 
\begin{enumerate}
	\item Modern categorical literature relies on strict notion of equality inherent to classical foundations to obtain various \emph{coherence theorems} and \emph{strictification results}, which are often an important step towards major categorical results \cite{MR0985657}. On the other side, the intensional nature of identities diminishes the efficacy of strictified definitions, prompting the exploration of alternative proof methods in the pursuit of formalizations.
	\item As explained in \cref{subsec:math}, (higher) categories are often considered up to suitably chosen equivalences. It is hence standard mathematical practice to extend results and properties to equivalent objects. However, a naive formalization of categories would fail to consider the data associated with such equivalences.  This would hence deprive us of the ability to transport along identities and result in a missed opportunity to employ a potent formalization technique when working with categories. 
	\item Many important examples of categories, such as the category of groups or rings, are defined in several equivalent ways. One significant component of modern categorical literature deftly uses these diverse characterizations to unveil new results and connections among established mathematical theories. However, imposing unified definitions diminishes our ability to delve into distinct equivalent, but not equal, facets of a mathematical theory.
\end{enumerate}

\subsection{The Univalence Maxim for Category Theory}
\label{subsec:univ}
Univalent foundations is a newly developed foundation for mathematics that relies on the univalence principle and, in contrast to classical foundations, offers a richer (\emph{intensional}) notion of identity. This stronger notion of equality has already resulted in a prominent role for univalent foundations in another branch of mathematics which also exhibits various notions of equivalences, namely \emph{topology}, resulting in \emph{synthetic algebraic topology} \cite{hottbook}. 

The main insight in this paper consists of demonstrating that defining categorical structures in this richer foundation can remedy the drawbacks outlined in \cref{subsec:math,subsec:cs}, via the following \emph{univalence maxim}: for every chosen equivalence of a given categorical structure, there exists a tailored definition whose notion of equality in univalent foundation precisely coincides with the chosen notion of equivalence. This means that in this setting we do not only have one notion of category or $2$-category, but rather one for each notion of equivalence.

Our concrete results, which due to their technical nature have been relegated to \cref{subsec:results}, consist of demonstrating manifestations of this maxim in the context of category theory, $2$-category theory and double category theory. Hence, we will focus here on how this maxim can overcome all the challenges outlined above. 
\begin{enumerate}
	\item We avoid any ambiguity by making the equivalence part of the definition of a categorical structure. 
	\item By severing the tie between equivalences and the structure of the categorical context, we can combine higher structures, such as natural transformations, and equivalences, such as isomorphisms of categories, based on our mathematical aims.
	\item Formalizing a higher categorical structure coupled with a designated equivalence obviates the need for strictification results, as instead the notion of equality has been adapted according to the equivalence. 
	\item Formalizing a categorical notion with a built-in notion of equivalence provides a direct mechanism for transporting along the specified equivalences, hence recovering the ability to use equivalent categories interchangeably. This, moreover, enables a seamless transition between different equivalent formalizations of one definition, recovering common mathematical practice.
    
    For example, with such principles at hand, one can simplify the characterization of invertible 2-cells and equivalences in complicated bicategories such as the bicategory of double categories \cite[Theorems 7.1 and 7.3]{DBLP:conf/cpp/WeideRAN24} and the bicategory of monads \cite[Propositions 3.6 and 3.7]{DBLP:conf/fscd/Weide23}.
\end{enumerate}

\subsection{Main Results}
\label{subsec:results}
In this paper we apply the univalence maxim stated above to the case of categories, $2$-categories and double categories. Starting on the technically less involved side, in \cref{sec:cat-uf}, we review the work in \cite{rezk_completion} and analyze two notions of categories in univalent foundations: \emph{setcategories}, which are invariant under isomorphisms, and \emph{univalent categories}, which are invariant under equivalences. After that, we review the work in \cite{DBLP:journals/mscs/AhrensFMVW21} and study four notions of $2$-categories in univalent foundations:
\begin{itemize}
	\item \emph{$2$-setcategories}, invariant under isomorphisms ;
	\item \emph{bisetcategories}, invariant under isomorphisms ;
	\item \emph{univalent $2$-categories}, invariant under equivalence of underlying categories, which are isomorphisms of $2$-cells ;
	\item \emph{univalent bicategories}, invariant under biequivalences.
\end{itemize}

The remainder of the paper focuses on the study of double categories. Double categories were originally introduced with the aim of better studying categories via \emph{formal category theory} \cite{MR0675339,MR0794752}, however, have now found a variety of applications in applied mathematics and computer science; see, for instance, its applications in systems theory \cite{courser2020open,Jaz_Myers_2021,Baez2022structuredversus}, programming languages theory \cite{DBLP:conf/lics/DagandM13,DBLP:conf/fossacs/NewL23}, rewriting \cite{DBLP:journals/jlap/BrownPHF23,DBLP:journals/jlap/BehrHK23,DBLP:conf/fscd/BehrMZ23}, databases \cite{DBLP:journals/corr/SchultzSVW16}, and lenses \cite{Clarke2023}. Double categorical structures exhibit a notoriously long and complicated list of structures and equivalences throughout the literature. This includes, but is not restricted to, \emph{strict double categories}, \emph{pseudo double categories}, \emph{double bicategories} with equivalences such as \emph{horizontal equivalences}, \emph{vertical biequivalences} and \emph{gregarious equivalences}. Applying our maxim to these double categorical notions permits us to clarify the relationship between the various structures and equivalences. 

More explicitly, 
\begin{enumerate}
	\item in \cref{sec:double-setcat} we introduce \emph{(pseudo) double setcategories} as (pseudo) double categories invariant under isomorphisms;
	\item in \cref{sec:univ-psdouble-cat} we analyze \emph{univalent pseudo double categories} invariant under vertical equivalences, reviewing work of \cite{DBLP:conf/cpp/WeideRAN24};
	\item we provide a further generalization to \emph{double bicategories} (\cref{sec:double-bicat}) and \emph{univalent double bicategories} invariant under gregarious equivalences (\cref{sec:univ-double-bicat}).
\end{enumerate} 

Finally, let us note that beyond clarifying the relations between double categorical characterizations and their equivalences, we also develop for the first time important aspects of the theory of double bicategories, such as companions (\cref{sec:comp-pair}) and computationally feasible methods to establish univalence (\cref{sec:horiz-inv}).

\subsection{A Diagrammatic Summary}
Given the technical nature of the definitions and their connections, we present here a diagrammatic summary of the major double categorical notions, how they relate to each other, and how they relate to different $2$-categorical structures. In the following diagram a dashed arrow represents an inclusion that only respects the categorical structure (i.e. only holds in classical setting), whereas a solid arrow indicates an inclusion that respects categorical structure and univalence conditions. Moreover, an arrow labeled $V$ denotes the underlying vertical $2$-category or bicategory (\cref{def:underlying-vert-strict-cat,rem:underlying-bicats}), whereas an arrow labeled $H$ is the underlying horizontal $2$-category or bicategory (\cref{def:underlying-hor-pseudo-double-cat,rem:underlying-bicats}).

\[
\begin{tikzcd}[column sep=-1cm]
	& \text{Strict Dbl Setcat (\ref{def:strict-double-setcat})} \arrow[dl, hookrightarrow] \arrow[dd, hookrightarrow] \arrow[rr, "HV" description] &[0.5cm] & \text{2-Setcat (\ref{def:settwocat})} \arrow[dd, hookrightarrow, dashed] \\
	\text{Ps Dbl Setcat (\ref{def:pseudo-double-set-cat})} \arrow[dr, hookrightarrow]  \arrow[dd, hookrightarrow, dashed]  \arrow[urrr, "V" description, bend right=10] \arrow[drrr, "H" description, bend left=10] & & &   \\
	& \text{Weak Dbl Setcat (\ref{rem:verity double bisetcategories})} \arrow[dd, hookrightarrow, dashed] \arrow[rr, "HV" description] & & \text{Bisetcat (\ref{def:settwocat})} \arrow[dd, hookrightarrow, dashed] \\
	\text{Univ Ps Dbl Cat (\ref{def:univalent-double-cat})} \arrow[dr, hookrightarrow, dashed] \arrow[ddr, bend right = 20, hookrightarrow, end anchor=north west] \arrow[rr, "V" description] \arrow[rrrd, dashed, "H" description] & & \text{Univ 2-Cat (\ref{def:univtwocat})} \arrow[ur, hookrightarrow, dashed] \arrow[dr, hookrightarrow, dashed] &  \\
	&  \text{Univ Weak Dbl Cat (\ref{def:univalent-weak-double-cat})} \arrow[d, hookrightarrow] \arrow[rr, "HV" description] &  &  \text{Univ Bicat (\ref{def:univbicat})} \\
	& \text{Univ Dbl Bicat (\ref{def:gregarious-univalence})} \arrow[urr, "HV" description, start anchor=north east] & & 
\end{tikzcd}
\]

This diagram only displays the interaction between the types of various double categorical notions, however we do anticipate a lift to inclusions of appropriately defined higher categories.

\subsection{Motivating Examples}
\label{subsec:exa}
Unsurprisingly a complicated diagram with this many definitions is motivated by a wide variety of examples in mathematics and computer science. Throughout the paper we hence introduce several classes of examples for each double categorical notion. Concretely, we analyze six important classes of examples that have found mathematical and computer scientific applications in formal category \cite{MR0675339,MR0794752}, rewriting \cite{DBLP:conf/fscd/BehrMZ23,DBLP:journals/jlap/BehrHK23}, open systems \cite{Baez2022structuredversus,MR4170469}, lenses \cite{DBLP:journals/pacmpl/BoisseauG18,DBLP:journals/corr/abs-2001-07488}, gradual type theory \cite{NewL19}, and ornaments \cite{DBLP:conf/lics/DagandM13}:
\begin{enumerate}
	\item \textbf{Squares:} \cref{exa:square-double-cat,exa:strict-square-double-cat,exa:symmetric-univalent-squares,exa:square-verity-double-bicat,exa:square-univ-double-bicat};
	\item \textbf{spans:} \cref{exa:span-pseudo-double-cat,exa:span-pseudo-double-setcat,exa:univalent-span};
	\item \textbf{structured cospans:} \cref{exa:struct-cospan-disp,exa:struct-cospan-disp-set,exa:struct-cospan-disp-univ}
	\item \textbf{profunctors:}  \cref{exa:profunctor-strict,exa:profunctor-strict-univalent,exa:enriched-profunctor-pseudo,exa:profunctor-verity-double-bicat,exa:enriched-profunctor-verity-double-bicat,exa:profunctor-nonstrict-not-univalent,exa:prof-double-bicat,exa:enriched-prof-double-bicat};
	\item \textbf{coreflections:} \cref{exa:coreflections,exa:coreflections-set,exa:coreflections-univ};
	\item \textbf{polynomial functors:} \cref{exa:polys,exa:polys-set,exa:polys-univ}.
\end{enumerate}
All of these classes are first defined in \cref{sec:def-pseudo-dbl-cat} and further analyzed and refined throughout the later sections. Beyond those six classes, we also consider examples coming from the \emph{mate calculus} \cite{MR0357542} in \cref{exa:mate,exa:mate-univ}.

\subsection{Related Work}
The study of categories in univalent foundations has a rich history. Indeed, a study of univalent categories was commenced in \cite{rezk_completion}, which was generalized to a study of univalence for $2$-categories and bicategories in \cite{DBLP:journals/mscs/AhrensFMVW21}, both of which are reviewed more carefully in \cref{sec:cat-uf}. Among the double categorical structures, univalent pseudo double categories and a univalence principle thereof is the focus of \cite{DBLP:conf/cpp/WeideRAN24}, which we review in \cref{sec:univ-psdouble-cat}. Finally, in \cite{up} the authors set up a general univalence principle, which in particular applies to univalent double bicategories, and is the motivation for \cref{sec:univ-double-bicat}.

Double categories acquired some attention from the formalization community, and several libraries on formalized mathematics contain double categories.
Murray, Pronk, and Szyld~\cite{doublecats_lean} worked towards defining double categories in the Lean proof assistant.%
\footnote{\url{https://github.com/leanprover-community/mathlib/pull/18204}}
In 1lab \cite{1lab}, internal categories are defined, and thus double categories are also defined as category objects in the category of setcategories.
Finally, in the library by Hu and Carette~\cite{DBLP:conf/cpp/HuC21} a definition of double category has been implemented\footnote{\url{https://github.com/agda/agda-categories/blob/36abe6bff98be027bd4fcc3306d6dac8b2140079/src/Categories/Double/Core.agda}}.
Each of these formalization only considers strict double categories, whereas we also consider weaker notions.
In addition, our formalization contains a study of various univalence conditions.

\section{(2-)Categories in Univalent Foundations}
\label{sec:cat-uf}
In this section we realize the vision of the univalence maxim (\cref{subsec:univ}) for categories and bicategories, based on work done in \cite{rezk_completion,DBLP:journals/mscs/AhrensFMVW21}. More precisely, we analyze two notions of equivalences for categories (isomorphisms, equivalences) and three notions of equivalences for $2$-categories and bicategories (isomorphisms, equivalences, biequivalences), and define for each one of the five equivalences precise categorical structures, such that their identities are the associated equivalences.

We will commence with categories. In classical mathematics a category is defined as a set of objects and a set of morphisms, depending on two objects, with a unital and associative composition of morphisms. In a type theoretical setting, we commence our characterization of categories in a similar manner as a type of objects and type of morphisms depending on two objects with a similarly characterized composition. However, in univalent foundations, only certain types exhibit a behavior similar to sets, meaning identities between two terms is a property rather than additional non-trivial structure. Hence, in order to get a category in univalent foundations, we need to specify that the morphisms form a set. A precise characterization of the data of a category in a type theoretical setting can be found in \cite[Definition 3.1]{rezk_completion}, where it is referred to as a \emph{pre-category}.

We now want to pursue the vision outlined in the introduction for categories. Categories have two important notions of equivalences:
\begin{enumerate}
	\item Isomorphisms of categories;
	\item equivalences of categories.
\end{enumerate}
We proceed to construct two notions of categories such that their identities correspond to these two equivalences. The first notion is straightforward, as it is the notion of sameness familiar from classical category theory. Indeed, by adding a further set-level assumption we get \emph{setcategories}, which are isomorphism invariant, as shown by Ahrens, Kapulkin, and Shulman \cite[Lemma 6.16]{rezk_completion}.

\begin{definition}[\coqident{CategoryTheory.Core.Setcategories}{setcategory}]
	\label{def:setcat}
	A category is said to be a \textbf{setcategory} if its type of objects is a set.
\end{definition}

\begin{remark}
	The notion of setcategory has already been introduced in \cite[Remark 3.19]{rezk_completion}, where it is called a \emph{strict category}.
\end{remark}

Obtaining an equivalence-invariant notion of category is more intricate and requires the notion of a \emph{univalent category}.

\begin{definition}[\coqident{CategoryTheory.Core.Univalence}{is_univalent}]
	\label{def:univ-cat}
	Let $\CC$ be a category. For all objects $x, y : \CC$ we have a map sending identities $p : x = y$ to isomorphisms $\idtoiso_{x, y}(p) : \iso{x}{y}$.
	A category is said to be \textbf{univalent} if for all objects $x, y : \CC$ the map $\idtoiso_{x, y}$ is an equivalence of types.
\end{definition}

The univalence condition now implies that an equality of categories is fully faithful and an equivalence of underlying groupoids, which is one characterization of an equivalence of categories, giving us the desired invariance under equivalences. See \cite[Theorem 6.17]{rezk_completion} for a more detailed argument. 

We can now use univalent categories to provide an alternative characterization of the invariance property of setcategories.

\begin{proposition}[\coqident{CategoryTheory.categories.CategoryOfSetCategories}{is_univalent_cat_of_setcategory}]
	\label{prop:univ-cat-setcat}
	The category of setcategories and functors is univalent, meaning that identities of setcategories correspond to isomorphisms. 
\end{proposition}

Unfortunately, we cannot repeat the argument and incorporate the equivalence invariance of univalent categories into the construction of a univalent category, as the type of equivalences between two univalent categories is generally too complex and does not form a set. However, we can in fact construct a univalent bicategory of univalent categories. 

\begin{proposition}[\coqident{Bicategories.Core.Examples.BicatOfUnivCats}{univalent_cat_is_univalent_2}]
	\label{prop:univ-bicat-univ-cat} 
	The bicategory of univalent categories, functors, and natural transformations is univalent, meaning that identities of univalent categories correspond to equivalences of categories (\cite[Theorem 6.17]{rezk_completion}) and that functor categories between univalent categories are univalent \cite[Theorem 4.5]{rezk_completion}.
\end{proposition}

See \cite[Proposition 3.19]{DBLP:journals/mscs/AhrensFMVW21} for a more detailed treatment of the bicategorical argument. 

Having exhausted all established notions of sameness of categories, we move on to bicategories, and $2$-categories, which are bicategories with identity associators and unitors. As bicategories not only have objects and $1$-morphisms, but also $2$-morphisms between $1$-morphisms, the number of relevant equivalences increases significantly. Instead of presenting an exhaustive list, here we focus in particular on three equivalences that are relevant in the coming sections: 
\begin{enumerate}
	\item isomorphisms of bicategories;
	\item equivalences of bicategories, meaning equivalences of the underlying $1$-categories that are isomorphisms of $2$-morphisms;
	\item biequivalences.
\end{enumerate}
We now construct for each type of equivalence a bicategory ($2$-category) such that their identities correspond to the chosen equivalence.

For the first kind of equivalence, taking \cref{def:setcat} as motivation we analogously impose appropriate set level restrictions to obtain isomorphism invariance. 
\begin{definition}[\coqident{Bicategories.Core.Strictness}{bisetcat}]
	\label{def:settwocat}
	A bicategory ($2$-category) is said to be a \textbf{bisetcategory} (\textbf{$2$-setcategory}) if its type of objects and of $1$-morphisms are sets.
\end{definition}

\begin{proposition}
	\label{prop:univ-cat-twosetcat}
	The two categories given by bisetcategories and functors, and by $2$-setcategories and functors are univalent, meaning that identities of $2$-setcategories and bisetcategories correspond to isomorphisms.
\end{proposition}

We now move on to the second kind of equivalence. Combining \cref{prop:univ-cat-setcat} and \cref{prop:univ-bicat-univ-cat}, we get the following result. 

\begin{definition}[\coqident{CategoryTheory.Core.TwoCategories}{univalent_two_cat}]
	\label{def:univtwocat}
	A $2$-category is said to be \textbf{univalent} if the underlying $1$-category is univalent and the $2$-morphisms form a set.
\end{definition}

\begin{proposition}
	\label{prop:loc_univ-cat-twosetcat}
	Identities of univalent $2$-categories correspond to equivalences of $2$-categories.
\end{proposition}

Finally, we want a $2$-categorical structure invariant under biequivalences, which is much more challenging. First, in light of \cref{prop:univ-bicat-univ-cat}, we need all hom categories to be univalent, which we call the \emph{local univalence} condition. Moreover, we also need a \emph{global univalence} condition, stating that identities in the type of objects correspond to equivalences in the $2$-category. In general this means that objects form a $2$-type and compositions of $1$-morphisms is generally not strictly associative or unital. Hence, this univalence condition only applies to bicategories. 

\begin{definition} \label{def:univbicat}
	A \emph{univalent bicategory} is a bicategory that is globally and locally univalent. 
\end{definition}
See \cite[Definition 3.1]{DBLP:journals/mscs/AhrensFMVW21} for a more explicit description of its definition. Finally, univalent bicategories do exhibit the anticipated invariance property; see \cite[Example 9.1]{up}.

\begin{proposition}
	\label{prop:univ-cat-univtwocat}
	Identities of univalent bicategories correspond to biequivalences. 
\end{proposition}

\section{Definition of Pseudo Double Categories} \label{sec:def-pseudo-dbl-cat}
In the previous section we reviewed (bi)categories and demonstrated how imposing additional conditions in univalent foundations results in (bi)categories up to a desired notion of sameness. For the remainder of this paper, we want to do a similar analysis for double categorical structures. As already explained in the introduction, the definition of a double category is far more intricate and involves more data, and so we have a much more diverse set of examples and equivalences. As a result, our analysis is both more non-trivial and mathematically valuable, as well as more applicable to the broader literature. 

In this transitional section we commence with a review of a general definition of pseudo double categories and explicate our examples introduced in \cref{subsec:exa}. Then, in the next two sections we show how imposing additional conditions result in the desired equivalences.

\begin{definition}[\coqident{Bicategories.DoubleCategories.Core.DoubleCats}{double_cat}]
	\label{def:double-cat-unfolded}
	A \textbf{pseudo double category} consists of
	\begin{enumerate}
		\item\label{double-cat:vertical-cat} a category $\CC$ called the \textbf{vertical category};
		\item\label{double-cat:hor-mor} for all objects $x : \CC$ and $y : \CC$, a type $\hor{x}{y}$ of \textbf{horizontal morphisms};
		\item\label{double-cat:hor-id} for every object $x : \CC$ a \textbf{horizontal identity} $\horid{x} : \hor{x}{x}$;
		\item\label{double-cat:hor-comp} for all horizontal morphisms $h : \hor{x}{y}$ and $k : \hor{y}{z}$, a \textbf{horizontal composition} $\horcomp{h}{k} : \hor{x}{z}$;
		\item\label{double-cat:squares} for all horizontal morphisms $h : \hor{x_1}{y_1}$ and $k : \hor{x_2}{y_2}$ and vertical morphisms $v_x : \hor{x_1}{x_2}$ and $v_y : \hor{y_1}{y_2}$, a set $\Dsquare{v_x}{v_y}{h}{k}$ of \textbf{squares};
		\item\label{double-cat:vertical-id-square} for all horizontal morphisms $h : \hor{x}{y}$, we have a \textbf{vertical identity} $\sqvid{h} : \Dsquare{\id{x}}{\id{y}}{h}{h}$;
		\item\label{double-cat:vertical-comp-square} for all squares $\tau_1 : \Dsquare{v_1}{w_1}{h}{k}$ and $\tau_2 : \Dsquare{v_2}{w_2}{k}{l}$,
		we have a \textbf{vertical composition}
		\[ \sqvcomp{\tau_1}{\tau_2} : \Dsquare{v_1 \cdot v_2}{w_1 \cdot w_2}{h}{l}; \]
		\item\label{double-cat:horizontal-id-square} for all $v : x \rightarrow y$,
		we have a \textbf{horizontal identity} \[ \sqhid{v} : \Dsquare{v}{v}{\horid{x}}{\horid{y}}; \]
		\item\label{double-cat:horizontal-comp-square} for all squares $\tau_1 : \Dsquare{v_1}{v_2}{h_1}{k_1}$ and $\tau_2 : \Dsquare{v_2}{v_3}{h_2}{k_2}$,
		we have a \textbf{horizontal composition}
		\[ \sqhcomp{\tau_1}{\tau_2} : \Dsquare{v_1}{v_3}{\horcomp{h_1}{h_2}}{\horcomp{k_1}{k_2}}; \]
		\item\label{double-cat:lunitor} for all $h : \hor{x}{y}$, we have a \textbf{left unitor} \[ \Dlunitor{h} : \Dsquare{\id{x}}{\id{y}}{\horcomp{\id{x}}{h}}{h}; \]
		\item\label{double-cat:runitor} for all $h : \hor{x}{y}$, we have a \textbf{right unitor} \[ \Drunitor{h} : \Dsquare{\id{x}}{\id{y}}{\horcomp{h}{\id{y}}}{h}; \]
		\item\label{double-cat:associator} for all $h_1 : \hor{w}{x}$, $h_2 : \hor{x}{y}$, and $h_3 : \hor{y}{z}$, we have an \textbf{associator} \[ \Dassociator{h_1}{h_2}{h_3} : \Dsquare{\id{w}}{\id{z}}{\horcomp{h_1}{(\horcomp{h_2}{h_3})}}{\horcomp{(\horcomp{h_1}{h_2})}{h_3}}. \]
	\end{enumerate}
	This data is required to satisfy several laws, stating, in particular,
	that horizontal identities and horizontal composition are functorial,
	and that the left unitor, right unitor, and associator are natural transformations.
	In addition, we have the \textbf{triangle} and \textbf{pentagon} law.
	Their description can be found in literature \cite{DBLP:conf/cpp/WeideRAN24}.
\end{definition}

 Upon observing the definition we notice that the horizontal and vertical morphisms are not defined symmetrically. Indeed, only the vertical morphisms have a strictly unital and associative composition. We can, however, slightly adjust the definition to obtain a stricter symmetric definition. 

\begin{definition}[\coqident{Bicategories.DoubleCategories.Core.StrictDoubleCats}{strict_double_cat}]
	\label{def:strict-double-cat}
	A \textbf{strict double category} is a pseudo double category in which horizontal morpshisms form a set and in which the unitors and associators of the horizontal composition are identities.
\end{definition}

We now proceed to give precise definitions of the four major classes of examples introduced in \cref{subsec:exa}: squares, spans, structured cospans and profunctors.
 \begin{example}[\coqident{Bicategories.DoubleCategories.Examples.SquareDoubleCat}{strict_square_double_cat}]
 	\label{exa:square-double-cat}
 	Let $\CC$ be a category.
 	Then we define a pseudo double category $\SquareS{\CC}$ as follows.
 	\begin{enumerate}
 		\item The objects are objects in $\CC$;
 		\item the horizontal and vertical morphisms are morphisms in $\CC$; and
 		\item the type of squares $\Dsquare{v_1}{v_2}{h_1}{h_2}$ is defined to be $h_1 \cdot v_2 = v_1 \cdot h_2$.
 	\end{enumerate}
 \end{example}
  
\begin{example}[\coqident{Bicategories.DoubleCategories.Examples.SpansDoubleCat}{spans_double_cat}]
	\label{exa:span-pseudo-double-cat}
	Let $\CC$ be a category with pullbacks.
	Then we define a pseudo double category $\SpanS{\CC}$ as follows.
	\begin{enumerate}
		\item The objects are objects in $\CC$;
		\item the vertical morphisms are morphisms in $\CC$;
		\item the horizontal morphisms are spans, which are diagrams of the form $\SpanDiag{x}{\varphi}{z}{\psi}{y}$; and
		\item given morphisms $f : x_1 \rightarrow x_2$ and $g : y_1 \rightarrow y_2$,
                  then a square with $f$ and $g$ as vertical sides and spans $\SpanDiag{x_1}{\varphi_1}{z_1}{\psi_1}{y_1}$ and $\SpanDiag{x_2}{\varphi_2}{z_2}{\psi_2}{y_2}$
                  as horizontal sides is a morphism $h : z_1 \rightarrow z_2$ such that the following
                  diagram commutes.
                  \[\begin{tikzcd}
                      {x_1} & {z_1} & {y_1} \\
                      {x_2} & {z_2} & {y_1}
                      \arrow["{\varphi_1}"', from=1-2, to=1-1]
                      \arrow["{\varphi_2}", from=2-2, to=2-1]
                      \arrow["{\psi_1}", from=1-2, to=1-3]
                      \arrow["{\psi_2}"', from=2-2, to=2-3]
                      \arrow["f"', from=1-1, to=2-1]
                      \arrow["g", from=1-3, to=2-3]
                      \arrow["h"{description}, from=1-2, to=2-2]
                    \end{tikzcd}\]
	\end{enumerate}
	Even if $\CC$ is a setcategory, $\SpanS{\CC}$ is generally only a pseudo double category.
        This is because composition of spans is given by pullbacks, which is only weakly unital and associative.
        Note that spans have been used in the study of rewriting systems \cite{DBLP:conf/fscd/BehrMZ23,DBLP:journals/jlap/BehrHK23}.
\end{example}

\begin{example}[\coqident{Bicategories.DoubleCategories.Examples.StructuredCospansDoubleCat}{structured_cospans_double_cat}]
\label{exa:struct-cospan-disp}
Suppose that we have a functor $L : \CC_1 \rightarrow \CC_2$.
We define the double category $\StructCospan{L}$ of structured cospans as follows.
\begin{enumerate}
  \item The objects are objects in $\CC_1$;
  \item the vertical morphisms are morphisms in $\CC_1$
  \item the horizontal morphisms are \textbf{structured cospans}, which are diagrams of the form $\Cospan{L(x)}{\varphi}{z}{\psi}{L(y)}$;
  \item given two structured cospans $\Cospan{L(x_1)}{\varphi_1}{z_1}{\psi_1}{L(y_1)}$ and $\Cospan{L(x_2)}{\varphi_2}{z_2}{\psi_2}{L(y_2)}$,
    and two morphisms $f : x_1 \rightarrow x_2$ and $g : y_1 \rightarrow y_2$,
    a square consists of a morphism $h : z_1 \rightarrow z_2$ such that the following diagram commutes
    \[\begin{tikzcd}
        {L(x_1)} & {z_1} & {L(y_1)} \\
        {L(x_2)} & {z_2} & {L(y_1)}
        \arrow["{\varphi_2}"', from=2-1, to=2-2]
        \arrow["{\psi_1}"', from=1-3, to=1-2]
        \arrow["{\psi_2}", from=2-3, to=2-2]
        \arrow["L(f)"', from=1-1, to=2-1]
        \arrow["L(g)", from=1-3, to=2-3]
        \arrow["h"{description}, from=1-2, to=2-2]
        \arrow["{\varphi_1}", from=1-1, to=1-2]
      \end{tikzcd}\]
\end{enumerate}
Note that structured cospans are used to study open systems \cite{Baez2022structuredversus,MR4170469}.  
\end{example}

\begin{example}[\coqident{Bicategories.DoubleCategories.Examples.ProfunctorDoubleCat}{strict_profunctor_double_cat}]
	\label{exa:profunctor-strict}
	We define the pseudo double category $\ProfDC$ as follows.
	\begin{enumerate}
		\item The objects are small setcategories;
		\item the vertical morphisms are functors;
		\item the horizontal morphisms are profunctors from a category $\CC$ to $\CD$, meaning functors of the form $\op{\CD} \times \CC \rightarrow \SET$;
		\item given profunctors $P : \profunctor{\CC_1}{\CD_1}$ and $Q : \profunctor{\CC_2}{\CD_2}$ and functors $F : \CC_1 \rightarrow \CC_2$ and $G : \CD_1 \rightarrow \CD_2$,
		we define squares $\Dsquare{F}{G}{P}{Q}$ to be natural transformations $\nattrans{P}{(F \times G) \cdot Q}$.
	\end{enumerate}
	Again, this pseudo double category will not be a strict double category, as profunctors do not compose strictly.
        Note that profunctors are important in the study of lenses \cite{DBLP:journals/pacmpl/BoisseauG18,DBLP:journals/corr/abs-2001-07488}.
\end{example}

Note that the composition of profunctors is defined as a colimit in the category of sets.
To guarantee that the desired colimit exists, we require that the involved categories are small.

\begin{example}[{Gradual type theory, \cite[Definition 5.2]{NewL19}}]
	\label{exa:coreflections}
	Let $\CC$ be a 2-category.
	We define a pseudo double category of coreflections, denoted $\Coreflect{\CC}$, as follows.
	\begin{enumerate}
		\item The objects are objects in $\CC$;
		\item vertical morphisms are morphisms in $\CC$;
		\item horizontal morphisms are adjunctions in $\CC$ whose unit is an equality (i.e., coreflections);
		\item squares are given by 2-cells in $\CC$.
	\end{enumerate}
\end{example}

\begin{example}[{Polynomial functors, \cite[Example 8]{DBLP:conf/lics/DagandM13}}]
\label{exa:polys}
Let $\CC$ be a locally Cartesian closed category.
We define a pseudo double category $\Polys{\CC}$ as follows.
\begin{enumerate}
  \item The objects are objects of $\CC$;
  \item the vertical morphisms are the morphisms of $\CC$;
  \item the horizontal morphisms from $I$ to $J$ are polynomial functors, i.e., diagrams of the form $I \leftarrow X \rightarrow Y \rightarrow J$ in $\CC$;
  \item squares are morphisms of polynomials.
\end{enumerate}
\end{example}

In \cite{DBLP:conf/lics/DagandM13} the authors define ornaments and descriptions, which are used to describe data types in Martin-L\"of Type Theory, and they define a pseudo double category of descriptions and ornaments that is equivalent to a subcategory of $\Polys{\SET}$ \cite[Proposition 2]{DBLP:conf/lics/DagandM13}, demonstrating the importance of polynomial functors in type theory. As part of their efforts, they develop several alternative descriptions of $\Polys{\SET}$, such as via slice categories and functors between them \cite[Example 8]{DBLP:conf/lics/DagandM13}.

Double categories are an intricate categorical structure who inherently include the data of several (2-)categories. In the next sections, while studying equivalences of double categories, we also want to analyze its compatibility with various underlying $2$-categories, necessitating precise definitions.

\begin{definition}[\coqident{Bicategories.DoubleCategories.Underlying.VerticalTwoCategory}{underlying_vert_two_cat}]
	\label{def:underlying-vert-strict-cat}
	Given a pseudo double category $\DC$, we define a strict 2-category $\VertCat{\DC}$, which we call the \textbf{underlying vertical 2-category}, as follows.
	\begin{enumerate}
		\item The objects are objects in $\DC$;
		\item the 1-cells are vertical morphisms in $\DC$; and
		\item the 2-cells are squares with horizontal identity sides.
	\end{enumerate}
\end{definition}

\begin{definition}[\coqident{Bicategories.DoubleCategories.Underlying.HorizontalBicategory}{horizontal_bicat}]
	\label{def:underlying-hor-pseudo-double-cat}
	Given a pseudo double category $\DC$, we define a bicategory $\HorCat{\DC}$, which we call the \textbf{underlying horizontal bicategory}, as follows.
	\begin{enumerate}
		\item The objects are objects in $\DC$;
		\item the 1-cells are horizontal morphisms in $\DC$; and
		\item the 2-cells are squares with vertical identity sides.
	\end{enumerate}
\end{definition}

\section{(Pseudo) Double Set-categories}
\label{sec:double-setcat}
Having defined strict and pseudo double categories, we now impose conditions to obtain an isomorphism invariant definition, following the vision of the univalence maxim (\cref{subsec:univ}). Concretely, we first construct isomorphism invariant notions of (pseudo) double categories and then study several classes of examples outlined in \cref{subsec:exa}. Following the discussions in \cref{sec:cat-uf}, and particularly \cref{def:setcat}, we need to impose a set level condition to obtain isomorphism invariance, motivating the following definitions.

\begin{definition}[\coqident{Bicategories.DoubleCategories.Core.StrictDoubleCats}{strict_double_cat}]
\label{def:strict-double-setcat}
A \textbf{strict double setcategory} is a strict double category whose objects form a set.
\end{definition}

\begin{definition}[\coqident{Bicategories.DoubleCategories.Core.PseudoDoubleSetCats}{pseudo_double_setcat}]
\label{def:pseudo-double-set-cat}
A \textbf{pseudo double setcategory} is a double category whose objects and horizontal morphisms form a set.
\end{definition}

Using similar ideas to the one used in \cref{prop:univ-cat-setcat}, we can now confirm their desired invariance.

\begin{theorem}[\coqident{Bicategories.DoubleCategories.Core.CatOfStrictDoubleCats}{univalent_cat_of_strict_double_cat}]
	\label{thm:univ-cat-of-strict-double-cat}
	The category of strict double setcategories, with objects being strict double setcategories and morphisms being strict double functors, is univalent.
        More concretely, identity of strict double setcategories corresponds with isomorphism.
\end{theorem}

\begin{theorem}[\coqident{Bicategories.DoubleCategories.Core.CatOfPseudoDoubleCats}{univalent_cat_of_pseudo_double_setcategory}]
	\label{thm:univ-cat-of-pseudo-double-cat} 
	The category of pseudo double setcategories, with objects being pseudo double setcategories and morphisms being strict double functors, is univalent.
        Concretely, this means that identity of pseudo double set categories corresponds with isomorphism.
\end{theorem}

We now review our motivating examples (\cref{subsec:exa}) in light of this invariance property.

\begin{example}[\coqident{Bicategories.DoubleCategories.Examples.SquareDoubleCat}{strict_square_double_cat}]
\label{exa:strict-square-double-cat}
Let $\CC$ be a setcategory, then $\SquareS{\CC}$ is also a strict double setcategory.
\end{example}

\begin{example}[\coqident{Bicategories.DoubleCategories.Examples.SpansDoubleCat}{spans_pseudo_double_setcat}]
  \label{exa:span-pseudo-double-setcat}
  If $\CC$ is a setcategory with pullbacks, then $\SpanS{\CC}$ is a pseudo double setcategory. 
\end{example}

\begin{example}[\coqident{Bicategories.DoubleCategories.Examples.StructuredCospansDoubleCat}{structured_cospans_pseudo_double_setcat}]
\label{exa:struct-cospan-disp-set}
If $\CC_1$ and $\CC_2$ are setcategories such that $\CC_2$ has pushouts,
then $\StructCospan{L}$ is a pseudo double setcategory.
\end{example}

\begin{example}
	\label{exa:coreflections-set}
	Let $\CC$ be a 2-setcategory.
	Then $\Coreflect{\CC}$ as defined in \cref{exa:coreflections} is a strict double setcategory.
\end{example}

\begin{example}
	\label{exa:polys-set}
	Let $\CC$ be a locally Cartesian closed setcategory. Then $\Polys{\CC}$, defined in \cref{exa:polys}, is a pseudo double setcategory.
\end{example}

The previous examples are very consistent with our intuition that an isomorphism invariant category should indeed result in isomorphism invariant double categories of squares, (co)spans, coreflections, or polynomials. Observe the class of examples of profunctors, \cref{exa:profunctor-strict}, is in fact not a pseudo double setcategory. See \cref{exa:profunctor-strict-univalent} for a more detailed discussion. Let us instead present one further example coming from $2$-category theory.

\begin{example}[\coqident{Bicategories.DoubleCategories.Examples.BiSetcatToDoubleCat}{bisetcat_to_pseudo_double_cat}]
	\label{exa:bicat-to-pseudo-double-setcat}
	Let $\CB$ be a bisetcategory.
	Then we define the pseudo setcategory $\BicatToDoubleCat{\CB}$ as follows.
	\begin{enumerate}
		\item The objects are objects in $\CB$;
		\item the vertical morphisms are morphisms in $\CB$;
		\item the horizontal morphisms from $x$ to $y$ are identities $x = y$; and
		\item the squares $\Dsquare{f}{g}{p}{q}$ are 2-cells
                  \[\cell{f \cdot \idtoiso(q)}{\idtoiso(p) \cdot g}.\]
	\end{enumerate}
\end{example}

We end this section by observing that the set condition is preserved by taking underlying bicategories. This results supports the intuitive fact that taking underlying $2$-categories preserves isomorphism invariance.

\begin{proposition}
  \label{prop:underlying-strict-cat}
  The following statements hold.
  \begin{enumerate}
    \item (\coqident{Bicategories.DoubleCategories.Underlying.VerticalTwoCategoryStrict}{strict_underlying_vertical_two_setcat})
      If $\DC$ is a strict double setcategory, then $\VertCat{\DC}$ is a $2$-setcategory. %
    \item (\coqident{Bicategories.DoubleCategories.Underlying.VerticalTwoCategory}{underlying_vert_two_setcat})
      If $\DC$ is a pseudo double setcategory, then $\VertCat{\DC}$ is a strict $2$-setcategory. %
    \item (\coqident{Bicategories.DoubleCategories.Underlying.HorizontalBicategory}{horizontal_setbicat})
      If $\DC$ is a pseudo double setcategory, then $\HorCat{\DC}$ is a bisetcategory. %
  \end{enumerate}
\end{proposition}

\section{Univalent Pseudo Double Categories}
\label{sec:univ-psdouble-cat}
In the previous section we focused on isomorphism invariance. In this section we continue realizing our univalence maxim (\cref{subsec:univ}), this time studying pseudo double categories invariant under \emph{vertical equivalences}, which are characterized by inducing equivalences on the underlying vertical category and for any two objects $x,y$ inducing equivalences on the category given by horizontal morphisms $\hor{x}{y}$ and squares with trivial vertical sides. Our analysis relies on work done in \cite{DBLP:conf/cpp/WeideRAN24}. We then end this section analyzing several important examples, following the pattern from \cref{subsec:exa}.

Following \cref{def:univ-cat}, we would expect a univalence condition for these two underlying categories. This, however, implies that the underlying horizontal category is neither univalent nor has a set of objects, meaning horizontal composition cannot be strict. Thus we have to work with pseudo double categories, resulting in the following definition.

\begin{definition}[\coqident{Bicategories.DoubleCategories.Core.UnivalentDoubleCats}{univalent_double_cat_to_double_cat}]
\label{def:univalent-double-cat}
A pseudo double category $\DC$ is said to be \textbf{univalent} if its underlying vertical category is univalent and if for all $x, y : \DC$ the category whose objects are horizontal morphisms $\hor{x}{y}$ and whose morphisms are squares with vertical identity sides, is univalent.
\end{definition}

Building on our insights in \cref{prop:univ-bicat-univ-cat}, we similarly construct a univalent bicategory of univalent pseudo double categories \cite{DBLP:conf/cpp/WeideRAN24}.

\begin{theorem}[\coqident{Bicategories.DoubleCategories.Core.BicatOfDoubleCats}{is_univalent_2_bicat_of_double_cats}]
\label{thm:bicat-univalent-double-cat}
The bicategory of univalent pseudo double categories with lax double functors is univalent.
Concretely, this means that identity of univalent pseudo double categories corresponds to adjoint equivalences.
\end{theorem}

Note that we use lax double functors in \cref{thm:bicat-univalent-double-cat} whereas we use strict double functors in \cref{thm:univ-cat-of-strict-double-cat,thm:univ-cat-of-pseudo-double-cat}.
As the univalence condition is motivated by vertical equivalences, it is not symmetric. For examples identities of objects only correspond to vertical isomorphisms, and identities of horizontal morphisms correspond to isomorphisms of squares (composed vertically).
However, some double categories satisfy a stronger univalence condition that is in fact symmetric.

\begin{definition}[\coqident{Bicategories.DoubleCategories.Core.SymmetricUnivalent}{symmetric_univalent}]
\label{def:symmetric-univalent}
A pseudo double category $\DC$ is said to be \textbf{symmetrically univalent} if the following conditions are satisfied.
\begin{enumerate}
  \item The horizontal morphisms form a set;
  \item $\DC$ is univalent;
  \item the category of objects and horizontal morphisms is univalent; and
  \item for all $x, y : \DC$ the category of vertical morphisms $\ver{x}{y}$ and squares with horizontal identity sides, is univalent.
\end{enumerate}
\end{definition}

Let us present a variety of examples of univalent and symmetrically univalent pseudo double categories, coming from our list in \cref{subsec:exa}.
Again we focus on our three classes of interest, namely squares, spans and profunctors \cite{DBLP:conf/cpp/WeideRAN24},
but we also provide additional examples.

\begin{example}[\coqident{Bicategories.DoubleCategories.Examples.SquareDoubleCat}{square_univalent_double_cat}]
\label{exa:symmetric-univalent-squares}
If $\CC$ is a univalent category, then $\SquareU{\CC}$ is a symmetrically univalent pseudo double category.
\end{example}

\begin{example}
\label{exa:coreflections-univ}
If $\CC$ is a univalent 2-category, then $\Coreflect{\CC}$ (\cref{exa:coreflections}) is a symmetrically univalent pseudo double category.
\end{example}

\begin{example}[\coqident{Bicategories.DoubleCategories.Examples.SpansDoubleCat}{spans_univalent_double_cat}]
\label{exa:univalent-span}
Let $\CC$ be a univalent category with pullbacks. Then the pseudo double category $\SpanS{\CC}$, is univalent, but not symmetrically univalent.
\end{example}

\begin{example}[\coqident{Bicategories.DoubleCategories.Examples.StructuredCospansDoubleCat}{structured_cospans_univalent_double_cat}]
\label{exa:struct-cospan-disp-univ}
If $\CC_1$ and $\CC_2$ are univalent categories such that $\CC_2$ has pushouts,
then $\StructCospan{L}$ is a univalent pseudo double category.
\end{example}

\begin{example}[\coqident{Bicategories.DoubleCategories.Examples.ProfunctorDoubleCat}{strict_profunctor_univalent_double_cat}]
\label{exa:profunctor-strict-univalent}
The pseudo double category of profunctors is univalent, as the category of setcategories is univalent (\cref{prop:univ-cat-setcat}), which also implies that it is not a pseudo double setcategory as there are non-trivial identities of objects given by isomorphisms of categories. Moreover, it is also not symmetrically univalent. Indeed, by \cite[Theorem 7.9.4]{borceux1994handbooki}, two categories are isomorphic in the category of categories and profunctors if their Cauchy completeness is isomorphic, which is a strictly weaker condition that being isomorphic.
\end{example}

In the later sections, we analyze enriched versions of double categories of profunctors benefiting from appropriately defined weak double categories (\cref{exa:enriched-profunctor-verity-double-bicat}). 
However, if we focus on categories enriched over a poset, which includes quantales and has found applications in automata theory \cite{MR1224222,MR1354679} and fuzzy logic \cite{MR4405034}, we do get stricter double categories. 

\begin{example}
\label{exa:enriched-profunctor-pseudo}
Suppose that $\CV$ is a complete and cocomplete symmetric monoidal category and suppose that $\CV$ is a poset.
Then we have a univalent pseudo double category whose objects are univalent categories enriched over $\CV$, vertical morphisms are enriched functors, horizontal morphisms are enriched profunctors,
and whose squares are given by enriched natural transformations.
\end{example}

\begin{remark}
	\label{exa:profunctor-nonstrict-not-univalent}
	Note we cannot construct a univalent pseudo double category given by univalent categories, functors and profunctors. Indeed, the type of univalent categories is a $2$-type as it includes all $1$-types, by \cite[Example 9.9.6]{hottbook} and \cite[Example 2.18]{DBLP:journals/mscs/AhrensFMVW21}, and hence cannot be the objects of a univalent category, which is at most a $1$-type (\cite[Lemma 3.8]{rezk_completion}). 
\end{remark}

\begin{example}
	\label{exa:polys-univ}
	Let $\CC$ be a locally Cartesian closed univalent category. Then $\Polys{\CC}$, defined in \cref{exa:polys}, is a univalent pseudo double category.
\end{example}

In the last section we observed that double setcategories (both pseudo and strict) preserve isomorphism invariance when taking underlying bicategories (\cref{prop:underlying-strict-cat}). We might hence anticipate similar results for univalent pseudo categories. Unfortunately we only have a partial result. Indeed, vertical equivalences are preserved by taking the underlying vertical $2$-category, confirmed by the following result.

\begin{proposition}[\coqident{Bicategories.DoubleCategories.Underlying.VerticalTwoCategory}{is_univalent_underlying_vert_two_cat}]
\label{prop:underlying-vertical-univalent-pseudo-double-cat}
If $\CD$ is a univalent pseudo double category, then $\VertCat{\CD}$ is a univalent $2$-category.
\end{proposition}

However, it is generally untrue that univalent pseudo double categories will induce globally univalent underlying horizontal bicategories. Indeed, the underlying horizontal bicategory of the pseudo double category $\ProfDC$ is given by small setcategories, profunctors and $2$-morphisms, and we already observed in \cref{exa:profunctor-strict-univalent} that profunctors can be an isomorphism without the underlying categories being isomorphic. However, not all hope is lost and we do recover a local univalence condition.

\begin{proposition}[\coqident{Bicategories.DoubleCategories.Underlying.HorizontalBicategory}{is_univalent_2_1_horizontal_bicat}]
\label{prop:underlying-horizontal-univalent-pseudo-double-cat}
If $\CD$ is a univalent pseudo double category, then $\HorCat{\CD}$ is locally univalent, but not necessarily globally univalent. 
\end{proposition}

This last result necessitates a double categorical notion that can accommodate biequivalences of bicategories, providing a first motivation for the coming sections.

\section{Motivating Verity Double Bicategories}
\label{sec:motivating-vdb}
Pseudo-double categories only exhibit non-strict compositions in the horizontal direction, so we anticipate a notion of a doubly weak double category with non-strict compositions in both directions. This is particularly relevant as pseudo double categories are unable to incorporate all relevant examples (\cref{exa:profunctor-nonstrict-not-univalent}) and have been unable to incorporate biequivalences as an invariant (\cref{prop:underlying-horizontal-univalent-pseudo-double-cat}). We could directly define such a generalized notion of double category, generalizing the constructors and coherences defined in \cref{def:double-cat-unfolded}, but this results in a new challenge. In a classical setting it would evidently follow that pseudo-categories fully faithfully embed in doubly weak double categories. However, in univalent foundations the situation becomes more challenging. Indeed, the strictness of vertical compositions in a univalent pseudo double category implies that identities in the type of objects correspond to vertical \textbf{isomorphisms}, whereas the weak vertical composition and symmetric nature of doubly weak double categories demands that identities corresponds to a pair of horizontal and vertical \textbf{equivalences} which interact well with each other, i.e. form a companion pair (\cref{sec:comp-pair}). Hence, a univalent pseudo double category will result in a non-univalent doubly weak double category. 

Given this conundrum, what we would need is a generalization of a double category, meaning a \emph{pre-double categorical structure}, whose notion of univalence can incorporate both univalent pseudo double categories and univalent doubly weak double categories. Taking the previous paragraph as motivation what such a structure should entail is a notion of $2$-morphisms divorced from the $2$-morphisms induced by squares (as defined in \cref{def:underlying-vert-strict-cat} and
\cref{def:underlying-hor-pseudo-double-cat}). By appropriately choosing the $2$-cells we can then restrict to both cases of interest: 
\begin{itemize}
	\item If we choose the vertical $2$-cells to be identities then equivalences in the vertical $2$-category would be isomorphism, recovering the identities of univalent pseudo-categories;
	\item if we choose the vertical (and horizontal) $2$-cells to coincide with $2$-cells induced by squares, then identities of objects correspond to equivalences we expect to see in doubly weak double categories. 
\end{itemize}

Hence, in order to further pursue our study of equivalences of double categories we first need to develop pre-double category theory, whose structure involves objects, horizontal morphisms, horizontal $2$-morphisms, vertical morphisms, vertical $2$-morphisms and squares, and their compositions along with an extensive list of coherences. Fortunately, a suitable candidate for such a notion has already been proposed by Verity in \cite{MR2844536}, where it is called a \emph{double bicategory}. Moreover its univalence principle has been analyzed in \cite{up}. Hence, for the remainder of the paper the aim is to study and formalize double bicategories, study its univalence principle, and establish an embedding from univalent pseudo categories to univalent double bicategories. 

\section{Verity Double Bicategories}
\label{sec:double-bicat}
Following the discussion of \cref{sec:motivating-vdb} we commence with the definition and formalization of double bicategories due to Verity \cite{MR2844536}, and also describe various examples, based on the framework presented in \cref{subsec:exa}. In \cref{sec:univ-double-bicat} we will then pursue its univalence principles. 

\begin{definition}[\coqident{Bicategories.DoubleCategories.DoubleBicat.VerityDoubleBicat}{verity_double_bicat}]
\label{def:verity-double-bicat}
A \textbf{Verity double bicategory} $\VB$ consists of
\begin{enumerate}
  \item a bicategory $\HorB[\VB]$ whose objects, 1-cells, and 2-cells are called \emph{horizontal};
  \item a bicategory $\VerB[\VB]$ with the same type of objects as $\HorB[\VB]$, and whose objects, 1-cells, and 2-cells are called \emph{vertical};
  \item for all objects $x_1, x_2, y_1, y_2 : \HorB[\VB]$, horizontal 1-cells $h_1 : \HorM{x_1}{x_2}$ and $h_2 : \HorM{y_1}{y_2}$, and vertical 1-cells $v_1 : \VerM{x_1}{y_1}$ and $v_2 : \VerM{x_2}{y_2}$ a set $\Dsquare{v_1}{v_2}{h_1}{h_2}$ of \textbf{squares};
  \item for all horizontal 1-cells $h : \HorM{x}{y}$ a square $\sqhid{h} : \Dsquare{\id{x}}{\id{y}}{h}{h}$;
  \item for all vertical 1-cells $v : \VerM{x}{y}$ a square $\sqvid{v} : \Dsquare{v}{v}{\id{x}}{\id{y}}$;
  \item\label{double-cat:vertical-comp-square-db} for all squares $s_1 : \Dsquare{v_1}{w_1}{h}{k}$ and $s_2 : \Dsquare{v_2}{w_2}{k}{l}$,
    we have a \textbf{vertical composition}
    \[ \sqvcomp{s_1}{s_2} : \Dsquare{v_1 \cdot v_2}{w_1 \cdot w_2}{h}{l}; \]
  \item\label{double-cat:horizontal-comp-square-db} for all squares $\tau_1 : \Dsquare{v_1}{v_2}{h_1}{k_1}$ and $\tau_2 : \Dsquare{v_2}{v_3}{h_2}{k_2}$,
    we have a \textbf{horizontal composition}
    \[ \sqhcomp{s_1}{s_2} : \Dsquare{v_1}{v_3}{\horcomp{h_1}{h_2}}{\horcomp{k_1}{k_2}}; \]
  \item given a vertical 2-cell $\tau : \VerC{v_1}{v_2}$ and a square $s : \Dsquare{v_2}{w}{h}{k}$, we have a \textbf{left whiskering} $\lwhisker{\tau}{s} : \Dsquare{v_1}{w}{h}{k}$;
  \item given a vertical 2-cell $\tau : \VerC{w_1}{w_2}$ and a square $s : \Dsquare{v}{w_1}{h}{k}$, we have a \textbf{right whiskering} $\rwhisker{\tau}{s} : \Dsquare{v}{w_2}{h}{k}$;
  \item given a horizontal 2-cell $\tau : \HorC{h_1}{h_2}$ and a square $s : \Dsquare{v}{w}{h_2}{k}$, we have a \textbf{up whiskering} $\uwhisker{\tau}{s} : \Dsquare{v}{w}{h_1}{k}$;
  \item given a horizontal 2-cell $\tau : \HorC{k_1}{k_2}$ and a square $s : \Dsquare{v}{w}{h}{k_1}$, we have a \textbf{down whiskering} $\dwhisker{\tau}{s} : \Dsquare{v}{w}{h}{k_2}$.
\end{enumerate}
In addition to the data explicated here, we have various laws governing their behavior. See the formalization or \cite[Definition 1.4.1]{MR2844536} for further details
\end{definition}

\begin{remark}
	\label{rem:underlying-bicats}
	Motivated by \cref{def:underlying-vert-strict-cat,def:underlying-hor-pseudo-double-cat}, for a given Verity double bicategory $\VB$, we call $\HorB[\VB]$ the \emph{underlying horizontal bicategory} and $\VerB[\VB]$ the \emph{underlying vertical bicategory}.
\end{remark}

Let us note that this definition does in fact satisfy all the desired conditions outlined in \cref{sec:motivating-vdb}. Indeed, we have independently defined horizontal and vertical $2$-cells and compositions of all $1$-morphisms are defined weakly, giving us a symmetric definition. Following our vision, we now define \emph{weak double categories} by adding an appropriate saturation condition identifying $2$-cells with certain squares.

\begin{definition}
	\label{def:saturated}
	Suppose we have a Verity double bicategory $\VB$.
	\begin{enumerate}
          \item (\coqident{Bicategories.DoubleCategories.DoubleBicat.CellsAndSquares}{horizontally_saturated})
            For all horizontal 1-cells $h_1, h_2 : \HorM{x}{y}$ we have a map $\CellToSquareH$ sending 2-cells $\tau : \HorC{h_1}{h_2}$ to the square $\uwhisker{\tau}{\sqhid{h_2}} : \Dsquare{\id{x}}{\id{y}}{h_1}{h_2}$.
            A Verity double bicategory is called \textbf{horizontally saturated} if the map $\CellToSquareH$ is an equivalence of types.
          \item (\coqident{Bicategories.DoubleCategories.DoubleBicat.CellsAndSquares}{vertically_saturated})
            Similarly, we have a map $\CellToSquareV$ sending vertical 2-cells $\tau : \VerC{v_1}{v_2}$ to the square $\lwhisker{\tau}{\sqvid{v_2}} : \Dsquare{v_1}{v_2}{\id{x}}{\id{y}}$,
            and we say that a Verity double bicategory is \textbf{vertically saturated} if the map $\CellToSquareV$ is an equivalence of types.
        \end{enumerate}
\end{definition}

\begin{definition}[\coqident{Bicategories.DoubleCategories.DoubleBicat.CellsAndSquares}{is_weak_double_cat}]
	\label{def:weak doublecat}
	A \emph{weak double category} is a horizontally and vertically saturated Verity double bicategory.
\end{definition}

Our definition of weak double categories by definition comes with an inclusion in Verity double bicategories. We now establish the second inclusion suggested in \cref{sec:motivating-vdb} and show that every pseudo double category gives rise to a Verity double bicategory.

\begin{example}[\coqident{Bicategories.DoubleCategories.Examples.DoubleCatToDoubleBicat}{double_cat_to_verity_double_bicat}]
	\label{exa:double-cat-to-verity-double-bicat}
	Suppose that we have a pseudo double category $\CC$.
	We define a Verity double bicategory $\DoubleCatToVerity{\CC}$ as follows.
	\begin{enumerate}
		\item The underlying horizontal bicategory is the discrete bicategory on the underlying vertical category of $\CC$;
		\item the underlying vertical bicategory is the underlying horizontal bicategory of $\CC$;
		\item squares are defined to be squares in $\CC$;
		\item the identity and composition operations are inherited from $\CC$;
		\item the left and right whiskering operations are defined using transport;
		\item the up and down whiskering operations are defined by taking a transport of a composition of squares.
	\end{enumerate}
	The Verity double bicategory $\DoubleCatToVerity{\CC}$ is vertically saturated, but not necessarily horizontally saturated.
\end{example}

Notice that the assignment of the horizontal $2$-cells in \cref{exa:double-cat-to-verity-double-bicat} is not unique, and our choice is motivated by the desire to realize our programme, meaning to guarantee that univalent pseudo-double categories give us univalent Verity double bicategories. See \cref{rem:choice-of-two-cells} for more details

We now proceed to look at several examples of Verity double bicategories, again following the pattern of examples introduced in \cref{subsec:exa}. By \cref{exa:double-cat-to-verity-double-bicat}, every pseudo double category results in a Verity double bicategory. So here we focus on examples giving us weak double categories, again motivated by our three classes of examples introduced in \cref{sec:def-pseudo-dbl-cat}, squares, profunctors and spans.  

\begin{example}[\coqident{Bicategories.DoubleCategories.Examples.SquareDoubleBicat}{square_verity_double_bicat}]
\label{exa:square-verity-double-bicat}
Let $\CB$ be a bicategory. We define a Verity double bicategory $\SquareV{\CB}$ as follows.
\begin{enumerate}
  \item The horizontal bicategory $\HorB[\SquareV{\CB}]$ is $\CB$;
  \item the vertical bicategory $\VerB[\SquareV{\CB}]$ is $\co{\CB}$;
  \item squares $\Dsquare{v}{w}{h}{k}$ are defined to be 2-cells $\cell{h \cdot w}{v \cdot k}$.
\end{enumerate}
We only show how the left whiskering operation is defined.
Suppose that we have 2-cells $\tau : \cell{h \cdot w}{v \cdot k}$ and $\theta : \cell{v}{v'}$,
then we define the left whiskering as the following composition of 2-cells
$\begin{tikzcd}
  {h \cdot w} & {v \cdot k} & {v' \cdot k.}
  \arrow["\tau", from=1-1, to=1-2, Rightarrow]
  \arrow["{\theta \vartriangleright k}", from=1-2, to=1-3, Rightarrow]
\end{tikzcd}$
Finally, note that $\SquareV{\CB}$ is both horizontally and vertically saturated, meaning it is a weak double category.
\end{example}
  
Note that in \cref{exa:square-verity-double-bicat}, the 2-cells in the vertical bicategory $\VerB[\SquareV{\CB}]$ are reversed.
This is necessary to get the right whiskering operations.

\begin{example}[\coqident{Bicategories.DoubleCategories.Examples.ProfunctorDoubleBicat}{univalent_profunctor_verity_double_bicat}]
	\label{exa:profunctor-verity-double-bicat}
	We define a Verity double bicategory $\ProfV$ as follows.
	\begin{enumerate}
		\item The underlying horizontal bicategory $\HorB[\ProfV]$ is the bicategory $\co{\bicatofcats}$ of small univalent categories;
		\item the objects of the underlying vertical $\VerB[\ProfV]$ are small univalent categories, 1-cells are profunctors, and 2-cells are natural transformations;
		\item the squares are defined in the same way as in \cref{exa:profunctor-strict}.
	\end{enumerate}
	The identity and composition operations for profunctors are defined in the same way as in \cref{exa:profunctor-strict}.
	Compared to \cref{exa:profunctor-strict}, we need to define whiskering operations and 2-cells of profunctors, and verify the axioms of Verity double bicategories.
	For the details of these constructions, we refer the reader to the formalization.
	Note that $\ProfV$ is both vertically and horizontally saturated, giving us a weak double category.
\end{example}

\begin{example}
\label{exa:enriched-profunctor-verity-double-bicat}
Let $\CV$ be a complete and cocomplete symmetric monoidal closed category.
We define a Verity double bicategory $\VProfV{\CV}$ as follows.
\begin{enumerate}
  \item The underlying horizontal bicategory $\HorB[\VProfV{\CV}]$ is the bicategory $\co{\bicatofEcats{\CV}}$ of small univalent enriched categories;
  \item the objects of the underlying vertical $\VerB[\VProfV{\CV}]$ are small univalent categories, 1-cells are enriched profunctors, and 2-cells are enriched natural transformations;
  \item the squares are defined in a similar way as in \cref{exa:profunctor-strict}.
\end{enumerate}
Identity, composition, and whiskering operations are defined similarly to \cref{exa:profunctor-verity-double-bicat}.
Since $\VProfV{\CV}$ is both vertically and horizontally saturated, it is a weak double category.
\end{example}

Note here we use small categories in \cref{exa:profunctor-verity-double-bicat,exa:enriched-profunctor-verity-double-bicat} to guarantee that the desired coends exist, analogous to \cref{exa:profunctor-strict}.

Unlike the previous sections we do not construct a weak double category of spans in a bicategory, as such a construction would necessitate an additional categorical layer. See \cite[Section 4]{morton2009cospans} for a possible approach. However, we do have one additional example motivated by \emph{mate calculus}.

\begin{example}[{Mate calculus, \cite[Proposition 2.2]{MR0357542}}]
	\label{exa:mate}
	Let $\CB$ be a bicategory.
	We define a Verity double bicategory $\LAdj{\CB}$ as follows.
	\begin{enumerate}
		\item The horizontal bicategory $\HorB[\LAdj{\CB}]$ is the bicategory whose objects are objects in $\CB$, 1-cells are left adjoints, and 2-cells are mate-pairs;
		\item the vertical bicategory $\VerB[\LAdj{\CB}]$ is $\co{\CB}$;
		\item given adjunctions $h : \adjunction{x}{y}$ and $k : \adjunction{x'}{y'}$, and 1-cells $v : x \rightarrow x'$ and $w ; y \rightarrow y$',
		squares $\Dsquare{v}{w}{h}{k}$ are defined to be 2-cells $h \cdot w \Rightarrow v \cdot k$.
	\end{enumerate}  
\end{example}

\begin{remark} \label{rem:verity double bisetcategories}
	Similar to \cref{sec:double-setcat}, we can impose a set condition on the the types of objects and morphisms to define \emph{Verity double bisetcategories} and \emph{weak double setcategories}, whose identities, analogous to \cref{thm:univ-cat-of-pseudo-double-cat,thm:univ-cat-of-strict-double-cat}, correspond to isomorphisms. However, similar to the classical setting (\cref{sec:motivating-vdb}), pseudo double setcategories fully faithfully embeds in weak double setcategories, hence obviating the need to generalize to any pre-double categorical notion, such as Verity double bicategories. 
	
	Given these similarities, we will hence directly proceed to the study of univalent Verity double bicategories and their relation to univalent pseudo categories and univalent weak double categories.
\end{remark}

\section{Companion Pairs}
\label{sec:comp-pair}
In the previous section we defined Verity double bicategories and showed that this general notion includes both pseudo double categories and weak double categories. Our major aim is to show that these assignments preserves univalence. Due to the inherently symmetric nature of Verity double bicategories, equivalences of objects are symmetric as well, meaning it needs to be given by a pair of horizontal and vertical equivalences that interact well with each other. In this section we provide a precise characterization of interaction between horizontal and vertical morphisms, via \emph{companion pairs.}

\begin{definition}[\coqident{Bicategories.DoubleCategories.DoubleBicat.CompanionPairs}{are_companions}]
\label{def:are-companions}
Suppose that we have a Verity double bicategory $\VB$ and a horizontal morphisms $h : \HorM{x}{y}$ and a vertical morphism $v : \VerM{x}{y}$.
Then we say that $h$ and $v$ form a \textbf{companion pair} if we have squares $\CompanionUnit : \Dsquare{v}{\id{y}}{h}{\id{y}}$ and $\CompanionCounit : \Dsquare{\id{x}}{v}{\id{x}}{h}$ such that the following squares $\rwhisker{\runitor}{\lwhisker{\linvunitor}{(\sqhcomp{\CompanionUnit}{\CompanionCounit})}}$ and $\dwhisker{\runitor}{\uwhisker{\linvunitor}{(\sqvcomp{\CompanionUnit}{\CompanionCounit})}}$ are identity squares.
We call $\CompanionUnit$ the \textbf{unit} and $\CompanionCounit$ the \textbf{counit}.
\end{definition}

Note an alternative notion of interaction between horizontal and vertical morphisms is given by \emph{conjoints}, which can be characterized as companion pairs in the horizontal dual of $\VB$. Historically speaking companion pairs arose in the study of double categories prevalent in \emph{formal category theory} \cite{DBLP:conf/fossacs/NewL23,MR0675339,MR0794752}.

Beyond the definition we also need several key properties of companion pairs that we confirm here. Concretely, we want to know that companion pairs include identities and are closed under composition. Moreover, companions, if they exist, are unique up to isomorphism. Finally, companions of equivalences are also equivalences.

\begin{example}[\coqident{Bicategories.DoubleCategories.DoubleBicat.CompanionPairs}{id_are_companions}]
\label{exa:identity-comp-companions}
Given an object $x$ in a Verity double bicategory $\VB$,
then the horizontal identity $\id{x}$ and vertical identity $\id{x}$ form a companion pair.

Suppose that $h_1 : \HorM{x}{y}$ and $v_1 : \VerM{x}{y}$ and that $h_2 : \HorM{y}{z}$ and $\VerM{y}{z}$ form companion pairs.
Then $h_1 \cdot h_2$ and $v_1 \cdot v_2$ also form a companion pair.
\end{example}

\begin{proposition}[\coqident{Bicategories.DoubleCategories.DoubleBicat.CompanionPairUnique}{isaprop_companion_pair}]
	\label{prop:unique-companions}
	Let $\VB$ be a Verity double bicategory such that $\HorB[\VB]$ and $\VerB[\VB]$ are locally univalent and such that $\VB$ is vertically saturated.
        For every horizontal 1-cell $h : \HorM{x}{y}$, we have that all vertical 1-cells $v, v' : \VerM{x}{y}$ that form a companion pair with $h$ are equal.
\end{proposition}

\begin{proposition}[\coqident{Bicategories.DoubleCategories.DoubleBicat.CompanionPairAdjEquiv}{companion_of_adjequiv}]
	\label{prop:companion-adj-equiv}
	Suppose that we have a Verity double bicategory $\VB$ such that $\VB$ is vertically saturated.
	Given a horizontal adjoint equivalence $\adjunction{l}{r}$ such that $l$ and $r$ have companion pairs $l'$ and $r'$ respectively,
	then we have a vertical adjoint equivalence given by $\adjunction{l'}{r'}$.
\end{proposition}

In many important examples of Verity double bicategories, every horizontal $1$-morphism has a companion, which, as we shall see in \cref{sec:horiz-inv}, is a key ingredient towards establishing its univalence condition. 
This holds for $\SquareV{\CB}$, $\ProfV$, and $\VProfV{\CV}$,
whereas if we have a pseudo double category $\CC$, then $\DoubleCatToVerity{\CC}$ has companion pairs if $\CC$ has.

\begin{example}[\coqident{Bicategories.DoubleCategories.Examples.SquareDoubleBicat}{all_companions_square_verity_double_bicat}]
\label{exa:square-companions}
Let $\CB$ be a bicategory.
Given a a 1-cell $f : x \rightarrow y$ in $\CB$, then $f$ and $f$ form a companion pair in $\SquareV{\CB}$.
\end{example}

\begin{example}[\coqident{Bicategories.DoubleCategories.Examples.ProfunctorDoubleBicat}{all_companions_univalent_profunctor_verity_double_bicat}]
\label{exa:profunctor-companions} 
Suppose that we have a functor $F : \CC_1 \rightarrow \CC_2$.
Note that $F$ gives rise to a profunctor $\reprL{F} : \profunctor{\CC_1}{\CC_2}$ that sends objects $x : \CC_1$ and $y : \CC_2$ to the set of morphisms $F(y) \rightarrow x$.
Then $F$ and $\reprL{F}$ form a companion pair.
\end{example}

\section{Univalent Double Bicategories}
\label{sec:univ-double-bicat}
In this section we use companion pairs introduced in \cref{sec:comp-pair} to present a univalence principle for double bicategories (\cref{sec:double-bicat}), further advancing our general maxim introduced in \cref{subsec:univ}. Given the amount data a Verity double bicategory involves, we split up the univalence condition into two parts. The first one is a local conditions imposed on the hom-categories in the underlying horizontal and vertical bicategories.

\begin{definition}[\coqident{Bicategories.DoubleCategories.DoubleBicat.LocalUnivalence}{locally_univalent_verity_double_bicat}]
\label{def:local-univalence}
A Verity double bicategory $\VB$ is said to be \textbf{locally univalent} if both $\HorB[\VB]$ and $\VerB[\VB]$ are locally univalent.
\end{definition}

The second univalence condition is global and focuses on the type of objects in Verity double bicategories. Here we use our newly gained understanding of companion pairs developed in \cref{sec:comp-pair}. 
\begin{definition}[\coqident{Bicategories.DoubleCategories.DoubleBicat.GregariousEquivalence}{gregarious_equivalence}]
	\label{def:gregarious-equiv}
	A \textbf{gregarious equivalence} from $x$ to $y$ in a Verity double bicategory $\VB$ consists of a horizontal adjoint equivalence $h : \HorM{x}{y}$, a vertical adjoint equivalence $v : \VerM{x}{y}$ such that $h$ and $v$ form a companion pair.
\end{definition}

\begin{example}[\coqident{Bicategories.DoubleCategories.DoubleBicat.GregariousEquivalence}{id_is_gregarious_equivalence}]
	\label{exa:id_gregarious-equiv}
	Given an object $x$ in a Verity double bicategory, then the horizontal identity $\id{x}$ and the vertical identity $\id{x}$ form a gregarious equivalence.
\end{example}

We now build on the insight in \cite{up} which defines \emph{gregarious univalence} using gregarious equivalences.

\begin{definition}[\coqident{Bicategories.DoubleCategories.DoubleBicat.GlobalUnivalence}{gregarious_univalent}]
\label{def:gregarious-univalence}
Given a Verity double bicategory and objects $x$ and $y$, we define the map $\idtogreg{x}{y}$ sending identities $x = y$ to gregarious equivalences
using path induction and the fact that the identity is a gregarious equivalence (\cref{exa:id_gregarious-equiv}).
A Verity double bicategory is said to be \textbf{gregarious univalent} if the map $\idtogreg{x}{y}$ is an equivalence of types for all $x$ and $y$.
\end{definition}

Note that in \cite[Example 9.3]{up} a univalence principle is proven for gregarious univalent Verity double bicategories.
More specifically, the identity type of such Verity double bicategories is equivalent to the type of gregarious equivalences between them.

Finally, \cref{def:weak doublecat} directly motivates the following definition.

\begin{definition}\label{def:univalent-weak-double-cat}
	A weak double category is \emph{univalent} if it is univalent as a double bicategory.
\end{definition}
This precisely coincides with the intuition presented in \cref{sec:motivating-vdb}, addressing the first inclusion. 

\section{Univalence and Weak Horizontal Invariance}
\label{sec:horiz-inv}
In this final section we tie up our discussion of univalent double categorical structures, by establishing the following two facts:
\begin{enumerate}
	\item There is a large class of univalent double bicategories coming from the classes of examples introduced in \cref{subsec:exa}, such as $\SquareV{\CB}$, $\ProfV$, and $\VProfV{\CV}$;
	\item every Verity double bicategory obtained from a univalent pseudo double category via the assignment \cref{exa:double-cat-to-verity-double-bicat} is univalent.
\end{enumerate} 
Given the complicated nature of gregarious univalence proving these directly is computationally challenging. We hence use a more conceptual approach. We define a notion of \emph{weakly horizontally invariant double bicategory} (\cref{def:horizontal-invariance}) and show in \cref{thm:gregarious-univalence-weq-horizontal-univalence} that in this case (with some minor conditions) gregarious univalence reduces to horizontal univalence. As a result, establishing gregarious univalence reduces to checking a horizontal invariance condition, as well as analyzing univalence of the underlying horizontal bicategory. We thus end this section with checking precisely these two properties for our cases of interest, namely the classes of examples introduced in \cref{subsec:exa} and univalent pseudo categories.

Let us start with the definition of weak horizontal invariance.

\begin{definition}[\coqident{Bicategories.DoubleCategories.DoubleBicat.CompanionPairs}{weakly_hor_invariant}]
\label{def:horizontal-invariance}
A Verity double bicategory is \textbf{weakly horizontally invariant} if every horizontal adjoint equivalence has a companion pair.
\end{definition}

\begin{example}[\coqident{Bicategories.DoubleCategories.DoubleBicat.CompanionPairs}{univalent_2_0_weakly_hor_invariant}]
\label{exa:global-univalence-horizontal-invariance}
Let $\VB$ be a Verity double bicategory such that $\HorB[\VB]$ is globally univalent.
Then $\VB$ is weakly horizontally invariant.
This is because the horizontal identity has a companion, and to construct companions for arbitrary adjoint equivalence, we use induction on adjoint equivalences.
\end{example}

Proving the main theorem requires the following proposition, which helps us characterize gregarious equivalences.
\begin{proposition}[\coqident{Bicategories.DoubleCategories.DoubleBicat.GregariousEquivalence}{hor_left_adjoint_equivalence_weq_gregarious_equivalence}]
\label{prop:gregarious-equivalence-invariance}
Let $\VB$ be a weakly horizontally invariant Verity double bicategory such that $\VB$ is vertically saturated and such that $\HorB[\VB]$ and $\VerB[\VB]$ are locally univalent.
Then a horizontal morphism $h : \HorM{x}{y}$ is an adjoint equivalence if and only if we have a vertical 1-cell $v : \VerM{x}{y}$ such that $h$ and $v$ are a gregarious equivalence.
\end{proposition}

\begin{theorem}[\coqident{Bicategories.DoubleCategories.DoubleBicat.GlobalUnivalence}{hor_globally_univalent_weq_gregarious_univalent}]
\label{thm:gregarious-univalence-weq-horizontal-univalence}
Let $\VB$ be a locally univalent, vertically saturated, and weakly horizontally invariant Verity double bicategory.
Then $\VB$ is gregarious univalent if and only if the bicategory $\HorB[\VB]$ is globally univalent.
\end{theorem}

We now apply \cref{thm:gregarious-univalence-weq-horizontal-univalence} to our cases of interest, by analyzing the classes of examples introduced in \cref{subsec:exa} as univalent double bicategories.

\begin{example}[\coqfile{Bicategories.DoubleCategories.Examples}{SquareDoubleBicat}]
\label{exa:square-univ-double-bicat}
Let $\CB$ be a univalent bicategory.
Note that both $\HorB[\SquareV{\CB}]$ and $\VerB[\SquareV{\CB}]$ are globally univalent.
Since the Verity double bicategory $\SquareV{\CB}$ is weakly horizontally invariant (\cref{exa:square-companions}), 
$\SquareV{\CB}$ is both locally univalent and gregarious univalent.
\end{example}

\begin{example}[\coqfile{Bicategories.DoubleCategories.Examples}{ProfunctorDoubleBicat}]
\label{exa:prof-double-bicat}
By \cite[Proposition 3.19]{DBLP:journals/mscs/AhrensFMVW21}, the bicategory of univalent categories is univalent.
Note that since identities of profunctors correspond to natural isomorphisms, the bicategory $\VerB[\ProfV]$ is locally univalent as well.
As such, the Verity double bicategory $\ProfV$ is both locally univalent and gregarious univalent, as it is weakly horizontally invariant, by \cref{exa:profunctor-companions}.
\end{example}

\begin{example}
\label{exa:enriched-prof-double-bicat}
	Following \cite[Theorem 2.6]{vanderweide2024enriched} $\VProfV{\CV}$ is horizontally univalent, and we can use a similar line of reasoning to \cref{exa:prof-double-bicat} to also conclude that $\VProfV{\CV}$ is weakly horizontally invariant. Hence, $\VProfV{\CV}$ is also gregarious univalent.	
\end{example}

\begin{example}
	\label{exa:mate-univ}
	Let $\CB$ be a univalent bicategory.
	We can show, analogous to \cref{exa:square-univ-double-bicat}, that $\LAdj{\CB}$ (\cref{exa:mate}) is locally univalent and gregarious univalent.
\end{example}

\begin{example}[\coqfile{Bicategories.DoubleCategories.Examples}{DoubleCatToDoubleBicat}]
\label{exa:univ-pseudo-double-cat}
Suppose that we have a univalent pseudo double category $\CC$.
Then the double bicategory $\DoubleCatToVerity{\CC}$ is gregarious univalent.
Indeed, weak horizontal invariance follows from \cref{exa:global-univalence-horizontal-invariance}.
The bicategory $\VerB[\DoubleCatToVerity{\CC}]$ is locally univalent,
and $\HorB[\DoubleCatToVerity{\CC}]$ is locally univalent by \cref{prop:underlying-horizontal-univalent-pseudo-double-cat}.
In addition, $\DoubleCatToVerity{\CC}$ is globally univalent because $\CC$ is univalent.
\end{example}

\begin{remark}
\label{rem:choice-of-two-cells}
Note, in \cref{exa:double-cat-to-verity-double-bicat} the choice of $2$-cells was not uniquely determined, and it was motivated by \cref{exa:univ-pseudo-double-cat}, as we explain. 
To obtain gregarious univalence, we need our choice of horizontal $2$-morphisms to be trivial so that identities are given by isomorphisms rather than more general equivalences. To give an explicit example, in the univalent pseudo-double category of setcategories, functors and profunctors \cref{exa:profunctor-strict}, the resulting Verity double bicategory is univalent because the $2$-cells in \cref{exa:double-cat-to-verity-double-bicat} are identities. If the vertical 2-cells would have been squares instead, both global and local univalence could fail. 

For example, in the case of the pseudo double category $\ProfDC$ we lose global univalence, as the existence natural transformations means that the identities in object need to correspond to equivalences of setcategories rather than isomorphisms thereof, and we lose local univalence, as the type of $1$-morphisms are sets and do not capture natural isomorphisms.  
\end{remark}

\begin{remark}
	The idea of using weak horizontal invariance to study equivalences is not new, see for example \cite{MR4602412}. From this perspective,  
	\cref{thm:gregarious-univalence-weq-horizontal-univalence} can be seen as a theoretical justification for the importance of weak horizontal invariance, by linking two different types of univalence.  
\end{remark}

\section{Formalization in UniMath}
The main results in this paper have been formalized using the Coq proof assistant
and the UniMath library.
The formalization uses 2-sided displayed categories, introduced by Van der Weide, Rasekh, Ahrens, and North in \cite{DBLP:conf/cpp/WeideRAN24}.

There are two differences between the formalization and the definitions presented in the paper.
First, the definition of strict double category in \cref{def:strict-double-cat} is presented as a special instance of pseudo double categories \cref{def:double-cat-unfolded} satisfying appropriate strictness conditions, whereas the formalization is given via an unfolded approach, along the lines of \cref{def:double-cat-unfolded}. 

Another deviation is that Verity double bicategories (\cref{def:verity-double-bicat}) are defined using two bicategories whose types of objects are equal.
However, we do not use that definition in the formalization: instead we use a trick to acquire two bicategories with the same collection of objects.
More specifically, we start with a bicategory $\HorB[\VB]$ and the data of a 2-sided displayed category to acquire vertical 1-cells, squares, vertical identity squares, and vertical composition.
Then we specify the necessary data and laws to acquire the bicategory $\VerB[\VB]$ of objects and vertical 1-cells,
and to do so, we reuse the definitions in \UniMath used to define bicategories.
These deviations between the formalization and the definition are motivated by the desire to obtain a more accessible and compact definition in the paper.

\section{Conclusion}
\label{sec:conclusion}
In this paper we presented a deep connection between equivalences of categorical structures and formalizations thereof in $\UniMath$, in the context of the univalence maxim for categorical structures (\cref{subsec:univ}). We discussed how, beyond being theoretically satisfying, it helps advance the study of category theory in mathematics and in computer science (\cref{subsec:results}), addressing various shortcomings in a classical setting (\cref{subsec:math,subsec:cs}). We back up our claim regarding the strength of this approach by analyzing this maxim in several categorical structures of interest, namely categories, $2$-categories and double categories. In \cref{sec:cat-uf} we focus on manifestations of this perspective in category theory and $2$-category theory. The remainder of our work focuses on a variety of double categorical structures. 

Concretely, we defined strict double setcategories (\cref{def:strict-double-setcat}) and pseudo double setcategories (\cref{def:pseudo-double-set-cat}) and showed that both are invariant under isomorphisms (\cref{thm:univ-cat-of-strict-double-cat,thm:univ-cat-of-pseudo-double-cat}). We then proceeded to define univalent pseudo categories (\cref{def:univalent-double-cat}) as pseudo double categorical notion invariant under horizontal equivalences (\cref{thm:bicat-univalent-double-cat}), building on work done by Van der Weide, Rasekh, Ahrens, and North \cite{DBLP:conf/cpp/WeideRAN24}.

Next, we take a major step and move from strict and pseudo double categories, which permit non-strict compositions at most in one direction, to Verity double bicategories (\cref{def:verity-double-bicat}), which should be thought of as a pre-double categorical notion that can incorporate both pseudo double categories (\cref{exa:double-cat-to-verity-double-bicat}) and weak double categories, which is a double category in which horizontal and vertical composition of morphisms is non-strict (\cref{def:weak doublecat}). We then move to our last central result and show that Verity double bicategories come with a notion of univalence, called gregarious univalence (\cref{def:gregarious-univalence}), which incorporates both univalent pseudo and weak double categories (\cref{exa:univ-pseudo-double-cat}). We obtain this result by establishing a new method of establishing gregarious univalence via weak horizontal invariance (\cref{thm:gregarious-univalence-weq-horizontal-univalence}). 

Throughout the whole paper, we analyzed many examples and how they manifest in specific circumstances, providing many distinguishing examples for all our definitions. A summary of all examples has already been presented in \cref{subsec:exa}.

Let us end with a note on possible future steps. For certain classes of (2-)categories and double categories, we obtain our univalence principle via the univalence of an appropriately defined (bi)category. See \cref{prop:univ-cat-setcat,prop:univ-bicat-univ-cat,thm:univ-cat-of-strict-double-cat,thm:univ-cat-of-pseudo-double-cat,thm:bicat-univalent-double-cat} for such results. Ideally we would like to to prove similar results for bicategories and Verity double bicategories, however, for that we would first need to develop \emph{tricategory theory} and a univalence principle thereof.

\begin{acks}
  We gratefully acknowledge the work by the Coq development team in providing the Coq proof assistant and surrounding infrastructure, as well as their support in keeping UniMath compatible with Coq.
  We are very grateful to Mike Shulman for answering our questions about profunctors.
  We would also like to thank Lyne Moser for many valuable discussions and explanations regarding double categories, their equivalences, and important references.
\end{acks}

\bibliographystyle{ACM-Reference-Format}
\bibliography{literature.bib}

\end{document}